\documentclass[12pt,leqno]{article}
\title{\ }
\author{\ }
\begin{document}
\newtheorem{theorem}{Theorem}[section]
\newtheorem{proposition}{Proposition}[section]
\newtheorem{lemma}{Lemma}[section]
\newtheorem{definition}{Definition}[section]
\newtheorem{corollaire}{Corollary}[section]
\newtheorem{remarque}{Remark}[section]
\def\Lim{\displaystyle\lim}
\def\Sup{\displaystyle\sup}
\def\Inf{\displaystyle\inf}
\def\Max{\displaystyle\max}
\def\Min{\displaystyle\min}
\def\Sum{\displaystyle\sum}
\def\Frac{\displaystyle\frac}
\def\Int{\displaystyle\int}
\def\n{|\kern -.05cm{|}\kern -.05cm{|}}
\def\Bigcap{\displaystyle\bigcap}
\def\E{{\cal E}}
\def\R{{\bf \hbox{\sc I\hskip -2pt R}}} 
\def\N{{\bf \hbox{\sc I\hskip -2pt N}}} 
\def\Z{{\bf Z}} 
\def\Q{{\bf \hbox{\sc I\hskip -7pt Q}}} 
\def\C{{\bf \hbox{\sc I\hskip -7pt C}}} 

\newcommand{\B}{\mathscr{B}}


{\begin{center} {\ } \vskip 2.1cm {\large\bf ON THE STABILITY OF DEGENERATE SCHR\"{O}DINGER EQUATION
UNDER BOUNDARY FRACTIONAL DAMPING\\ } \
\\ \ \\
\begin{tabular}{c}
\sc  Fatiha Chouaou and Abbes Benaissa,
\\ \\
\small  Laboratory of Analysis and Control of PDEs, \qquad \\
\small Djillali Liabes University,\qquad \\
\small P. O. Box 89, Sidi Bel Abbes 22000, ALGERIA.\qquad \\
\small   drchouaou@gmail.com\\
\small benaissa$_{-}$abbes@yahoo.com \\
\end{tabular}
\end{center}
\
\\
\begin{abstract}
{\small In this paper we study the well-posedness and stability of degenerate Schr\"{o}dinger equation with a
fractional boundary damping. First, we establish the well-posedness of the degenerate problem
$\psi_t(x,t)-\imath(\tau(x) \psi_x(x,t))_x=0, \hbox{ with } x \in (0,1)$, controlled by Dirichlet-Neumann conditions.
Then, exponential  and polynomial decay rate of the solution are established using multiplier method.}

\end{abstract}


{\it ${\cal AMS}$ Classification:} 35B40, 35Q41, 35P15.
\begin{sloppypar}

\section{Introduction}
In this paper, we study the well-posedness and stability of degenerate Schr\"{o}dinger equation under a fractional
boundary damping represented by the following system:
\begin{equation}\label{P2'}
\quad \left \{
\matrix{
\psi_t(x,t)-\imath(\tau(x) \psi_x(x,t))_x=0, \hfill &  \hbox{ in } (0,1)\times(0,+\infty),  \hfill &\cr
\left\{\matrix{\psi(0, t)=0 \hfill & \hbox{ if } 0\leq m_{\tau}< 1 \hfill & \cr
(\tau(x) \psi_{x})(0,t)=0 \hfill & \hbox{ if } 1\leq m_{\tau}< 2 \hfill & \cr}\right.\hfill  & \hbox{ in } (0, +\infty) ,\hfill&\cr
(\tau\psi_x)(1,t)-\imath\rho \partial^{\tilde{\alpha}, \wp}\psi(1,t)=0 \hfill & \hbox{ in } (0,+\infty),\hfill & \cr
\psi(x,0)=\psi_0(x) \hfill &    \hbox{ on } (0,1),\hfill & \cr}
\right.
\end{equation}
where $\rho>0$ and $\tau \in {\cal C}([0,1])\cap {\cal C}^1([0,1])$ a function  satisfying the following assumptions:
$$
\quad \left \{
\begin{array}{ll}
\tau(x)>0, \ \forall \ x \in ]0,1], \ \tau(0)=0,\\\\
m_{\tau}=\Sup_{x<0\leq 1} \Frac{x|\tau'(x)|}{\tau(x)}<2,\\\\
\tau \in {\cal C}^{[m_{\tau}]}([0,1]),
\end{array}
\right.
$$
where $[.]$ stands for the integer part.
$\partial^{\tilde{\alpha}, \wp}$ represents the exponential fractional integral of order $0<\tilde{\alpha}\leq 1$
(see {\bf\cite{choi}}), which is given by
$$
\partial^{\tilde{\alpha}, \wp}w({t})=\left\{\matrix{w(t)\hfill &\hbox{ for } \tilde{\alpha}=1,\ \wp\geq 0,\hfill \cr
\Frac{1}{\Gamma(1-\tilde{\alpha})}\Int_{0}^{{t}} ({t}-s)^{-\tilde{\alpha}}e^{-{\wp}({t}-s)}w(s)\, ds,\hfill & \hbox{ for } 0<\tilde{\alpha}<1,\ \wp\geq 0.\hfill \cr}\right.
\label{cap}
$$
We analyse the two following cases:
\begin{itemize}
\item[] When $m_{\tau}\in [0, 1[$, then the problem
is named weakly degenerate at $0$ and the usual boundary condition related to $(\ref{P2'})_1$ is
the Dirichlet boundary condition $\psi(0) = 0$.

\item[] When $m_{\tau}> 1$, then the problem
is named strongly degenerate at $0$ and the usual boundary condition related to $(\ref{P2'})_1$ is
the Neumann boundary condition $(\tau\psi_x)(0) = 0$.
\end{itemize}

Equation $(\ref{P2'})$ has diverse applications in the field of nonlinear
optics, plasma physics, fluid mechanics and quantum gases.
Taking into account such a type of damping is not only of mathematical interests but it is also
related to concrete real-life applications.

G. Fibich {\bf\cite{fibic}} acclaimed that the damping (absorption) term plays an
important effect in the physical model and it is better to not be
ignored. This inspired us to consider $(\ref{P2'})$.

Many authors addressed the Schr\"{o}dinger equation like Machtyngier and Zuazua {\bf\cite{MZ}},
where they studied the stabilization of the Schr\"{o}dinger equation by multiplier techniques.
they obtained an exponential decay under certain geometric conditions in $\Omega$ an open bounded domain in $\R^n$.

Liu and Wang {\bf\cite{LW}} applied the backsteping method to study the stabilization problem of an antistable
Schr\"{o}dinger equation by boundary feedback with only displacement observation.
Firstly, they obtained an exponential stability in $H^{1}(\Omega)$
with Newmann boundary damping, then they studied the stability in $L^{2}(\Omega)$ where the damping term is
in a neighborhood of board of $\Omega$.

Cui, Liu and Xu {\bf\cite{cuili}} traited the stabilization of Schr\"{o}dinger equation with
constrained boundary damping. These stability results are similar to those obtained in {\bf\cite{KZ}} and {\bf\cite{Z}}.

Nicaise and Rebiai {\bf\cite{nire}} obtained the exponential
stabilization of a Schr\"{o}dinger equation with boundary or internal difference-type
damping.
Recently, Chen, Xie and Xu {\bf\cite{chxi}} applied the approach comes from
Feng et al. {\bf\cite{fexu}} to the multi-dimensional Schr\"{o}dinger equation with interior pure delay
control and obtained the exponential stabilization. Moreover, Li, Chen and Xie {\bf\cite{liche}}
used this method to consider the stabilization problem for the Schr\"{o}dinger equation
with an input time delay, and obtained the exponential stability with arbitrary decay
rate of the closed-loop system.

This paper is organized as follows. In section 2, we give some preliminaries. In section 3,
the well-posedness results of the system $(\ref{P2'})$ are given using semigroup theory. In section 4,
we prove exponential and polynomial decay using classical multiplier method. In section 5, we show that the obtained results are optimal.

\section{Preliminaries}

\quad In this section, we introduce notations, definitions and
propositions that will be used later. First we introduce some  weighted Sobolev spaces:
$$
H^1_{\tau}(0,1)=\left\{\psi \in L^2(0,1), \psi\hbox{ is locally absolutely continuous in }(0,1], \ \sqrt{\tau(x)}\psi_x\in L^2(0,1)\right\},
$$
$$
H^1_{0,\tau}(0,1)=\left\{\psi \in H^1_{\alpha}(0,1), \ \psi(0)=0\right\},
$$
and
$$
H^2_{\tau}(0,1)=\left\{\psi \in L^2(0,1), \psi\in H^1_{\tau}(0,1), \ \tau(x)\psi_x\in H^1(0,1)\right\},
$$
where $H^1(0,1)$ represent the classical Sobolev space.\\

We can set the following norm
$$
|\psi|_{H^1_{0, \tau}(0,1)}=\left(\int_0^1\tau(x)|\psi_x(x)^2|dx\right)^{1/2}, \quad  \forall \ \psi \in H^1_{0, \tau}(0,1),
$$
which is an equivalent norm on the closed subspace $H^1_{0, \tau}(0,1)$ to the norm of $H^1_{\tau}(0,1)$.

In order to express the boundary conditions of the first component of the solution of
(\ref{P2'}) in the functional setting, we define the space $W_{\tau}^1(0, 1)$ depending on the value of $m_{\tau}$,
as follows:
$$
W_{\tau}^{1}(0,1)=\left\{\matrix{H^1_{0,\tau}(0,1) \hfill & \hbox{ if }0\leq m_\tau< 1,\hfill &\cr
H_{\tau}^{1}(0,1)\hfill & \hbox{ if }1\leq m_\tau< 2.\hfill &\cr}\right.
$$

\subsection{Augmented model}
In this section we reformulate $(\ref{P2'})$ into an augmented system. For that, we need the following proposition.
\begin{proposition}[see {\bf\cite{mbod}}]
Let $\eta$ be the function:
\begin{equation}
\eta(\xi)=|\xi|^{(2\tilde{\alpha}-1)/2},\quad -\infty<\xi<+\infty,\
0<\tilde{\alpha}<1. \label{e1}
\end{equation}
Then the relationship between the `input' U and the `output' O of
the system
\begin{equation}
\partial_t\vartheta(\xi, t)+(\xi^{2}+\wp)\vartheta(\xi, t) -U(t)\eta(\xi)=0,\quad -\infty<\xi<+\infty,\wp\geq 0, t> 0,
\label{e6}
\end{equation}
\begin{equation}
\vartheta(\xi, 0)=0, \label{e7}
\end{equation}
\begin{equation}
O(t)=(\pi)^{-1}\sin(\tilde{\alpha}\pi)\Int_{-\infty}^{+\infty}\eta(\xi)\vartheta(\xi, t)\, d\xi,
\label{e8}
\end{equation}
where $U\in C^0([0,+\infty))$, is given by
\begin{equation}
O=I^{1-\tilde{\alpha}, \wp}U, \label{e9}
\end{equation}
where
$$
[I^{\alpha,
\wp}f](t)=\Frac{1}{\Gamma(\tilde{\alpha})}\Int_{0}^{t}(t-s)^{\tilde{\alpha}-1}e^{-\wp(t-s)}f(s)\, ds.
$$
\label{th2}
\end{proposition}
\begin{lemma}[see {\bf\cite{achour.1}}]
If $\lambda\in D_{\wp}=\C\backslash]-\infty, -\wp]$ then
$$
\Int_{-\infty}^{+\infty}\Frac{\eta^2(\xi)}{\lambda+\wp+\xi^2}\, d\xi
=\Frac{\pi}{\sin\tilde{\alpha}\pi}(\lambda+\wp)^{\tilde{\alpha}-1}.
$$
\label{achour}
\end{lemma}
Using now Proposition \ref{th2} and relation (\ref{e9}), system $(\ref{P2'})$ may be rework into the following
augmented system
$$
\left\{\matrix{\psi_t(x,t)-\imath(\tau(x) \psi_x(x,t))_x=0,\hfill  & \hfill &\cr
\vartheta_t(\xi, t)+(\xi^{2}+\wp)\vartheta(\xi, t) -\psi(1, t)\eta(\xi)=0,\hfill & -\infty <\xi<+\infty,\ \ \hfill t>0,\hfill&\cr
\left\{\matrix{\psi(0, t)=0 \hfill & \hbox{ if } 0\leq m_{\tau}< 1 \hfill & \cr
(\tau(x)\psi_{x})(0,t)=0 \hfill & \hbox{ if } 1\leq m_{\tau}< 2 \hfill & \cr}\right.\hfill  & \hbox{ in } (0, +\infty) ,\hfill&\cr
(\tau\psi_x)(1,t)-\imath\zeta \Int_{-\infty}^{+\infty}\eta(\xi)\vartheta(\xi, t)\, d\xi=0 \hfill & \hbox{ in } (0,+\infty),\hfill & \cr
\psi(x,0)=\psi_0(x), \vartheta(\xi, 0)=0,\hfill &\cr }\right.
\leqno{(P)}
$$
where $\zeta=\rho(\pi)^{-1}\sin(\tilde{\alpha}\pi)$.

Furthermore, define the total energy for the system $(P)$ as
\begin{equation}\label{27'}
  E(t)=\frac{1}{2}\int_{0}^{1}|\psi(x,t)|^2\, dx+\frac{\zeta}{2}\Int_{-\infty}^{+\infty}|\vartheta(\xi, t)|^2\, d\xi.
\end{equation}
\begin{lemma}
Let $(\psi, \vartheta)$ be a regular solution of the problem $(P)$.
Then, the energy functional defined by (\ref{27'}) satisfies
\begin{equation}
E'(t)=-\zeta \Int_{-\infty}^{+\infty}(\xi^{2}+\wp)|\vartheta(\xi,t)|^2d\xi\leq 0.
\label{e11}
\end{equation}
\label{lem1}
\end{lemma}
{\bf Proof.}
Multiplying the first equation of $(P)$ by $\overline{\psi}$, integrating over $(0,1)$, applying integration by parts and
using boundary conditions we obtain
$$
\int_{0}^{1}\psi_t\overline{\psi} \ dx =- \zeta \Int_{\R}\eta(\xi)  \vartheta(\xi, t)\, d\xi\overline{\psi}(1, t)
-\imath \int_{0}^{1}\tau(x) |\psi_x|^2 \ dx,
\leqno{(s1)}
$$
Multiplying the second equation in $(P)$ by $\zeta \overline{\vartheta}$ and integrating over $(-\infty,+\infty)$, we get:
$$
\ \ \ \ \Frac{\zeta}{2}\Frac{d}{dt} \|\vartheta\|^{2}_{L^2(\R)}
+ \zeta \Int_{\R}(\xi^{2}+\eta)|\vartheta(\xi,t)|^2\, d\xi
- \zeta \Re\Int_{\R}\eta(\xi)  \ \overline{\vartheta}(\xi,t)d\xi \psi(1, t) = 0.
\leqno{(s2)}
$$
Taking the sum of the real part of $(s1)$ and $(s2)$, we get the result.
\section{Well-posedness}
The goal of this section is to prove the well-posedness results of $(\ref{P2'})$ using a semigroup approach and the
Lumer-Philips' Theorem.\\

We consider system $(P)$ in the natural state space
$$
{\cal H}=L^{2}(0,1)\times L^2(\R)
$$
equipped with the inner product
$$
\langle U, \widetilde{U} \rangle_{\cal H}=\int_{0}^{1}\psi(x)\overline{\widetilde{\psi}}(x)dx+\zeta\Int_{-\infty}^{+\infty}\vartheta(\xi)\overline{\widetilde{\vartheta}}(\xi)\, d\xi
$$
for all $U, \widetilde{U} \in {\mathcal{H}}$ with $U=(\psi, \vartheta)^{T}$ and $\widetilde{U}=(\widetilde{\psi}, \widetilde{\vartheta})^{T}$.
Define the system operator of $(P)$ as follows:
$$
{\cal A}U=\pmatrix{\imath (\tau(x) \psi_x)_x  \cr
-(\xi^2+\wp)\vartheta+\eta(\xi)\psi(1) \cr}
$$
with domain
$$
D({\mathcal{A}})=\left\{\matrix{ (\psi, \vartheta)\in {\cal H}: \psi\in H_{\tau}^{2}(0, 1)\cap W_{\tau}^1(0, 1), \
(\tau(x) \psi_x)(1)=\imath\zeta\Int_{-\infty}^{\infty}\eta(\xi)\vartheta(\xi)\, d\xi & \cr
-(\xi^2+\wp)\vartheta+\eta(\xi)\psi(1)\in L^2(\R), |\xi|\vartheta\in L^2(\R)& \cr}\right\}.
$$
Then, $(P)$ can be written as an evolution equation in ${\mathcal{H}}$
\begin{equation}\label{10'}
 \quad \left \{
\begin{array}{ll}
U_t={\cal A}U,\\
U(0)=U_0.
\end{array}
\right.
\end{equation}
The following theorem gives the existence and uniqueness results of the
solution of the problem (\ref{10'}).
\begin{theorem} Let $U_0 \in \mathcal{H}$, then there exists a
unique solution $U\in {{\cal C}}^0({\R}_+,{\cal H})$.Moreover, if $U_0 \in D({\cal A})$, then
$U\in {{\cal C}}^0({\R}_+,D({\cal A}))\cap {{\cal C}}^1(\R_+,{\cal H})$.
\label{1th}
\end{theorem}
{\bf Proof}\\
We use the semigroup approach. In what follows, we prove that ${\cal A}$ is dissipative. For any $U\in D({\cal A})$
and using (\ref{27'}), (\ref{e11}) and the fact that
\begin{equation}
E(t)=\Frac{1}{2}\|U\|_{\cal H}^{2},
\label{e16}
\end{equation}
we have
\begin{equation}
\Re\langle {\cal A} U,  U \rangle_{\cal H}=-\zeta \Int_{-\infty}^{+\infty}(\xi^{2}+\wp)|\vartheta(\xi)|^2d\xi\leq 0.
\label{e17}
\end{equation}
Hence, ${\cal A}$ is dissipative. Next, we prove that the operator $\lambda I-{\cal A}$ is surjective for $\lambda > 0$.
Given $F\in {\cal H}$, we prove that there exists $U\in D({\cal A})$ satisfying
\begin{equation}
(\lambda I-{\cal A})U=F.
\label{e18z}
\end{equation}
Equation (\ref{e18z}) is equivalent to
\begin{equation}
\left\{\matrix{\lambda \psi-\imath (\tau(x) \psi_x)_x=f_1,  \hfill & \cr
\lambda \vartheta+(\xi^{2}+\wp)\vartheta-\psi(1)\eta(\xi)=f_2.\hfill & \cr}\right.
\label{e18}
\end{equation}
By $(\ref{e18})_2$ we can find $\vartheta$ as
\begin{equation}
\vartheta(\xi)=\Frac{\psi(1)\eta(\xi)+f_2(\xi)}{\lambda +\xi^{2}+\wp}.
\label{m18}
\end{equation}
Problem (\ref{e18}) can reformulate as follows
\begin{equation}
\Int_{0}^{1}(\lambda \psi-\imath (\tau(x) \psi_x)_x) {\overline w}\, dx=\Int_{0}^{1} f_1{\overline w}\, dx \quad \forall w\in W_{\tau}^1(0, 1).
\label{eq55}
\end{equation}
Integrating by parts, we obtain
$$
\Int_{0}^{1}\lambda \psi {\overline w}\, dx+\rho (\lambda+\wp)^{\tilde{\alpha}-1}\psi(1)\overline{w}(1)
+\imath\Int_{0}^{1} \tau(x) \psi_x {\overline w}_x\, dx=\Int_{0}^{1} f_1{\overline w}\, dx-\zeta\Int_{-\infty}^{+\infty}\Frac{\eta(\xi)f_2(\xi)}{\lambda +\xi^{2}+\wp}\, d\xi \overline{w}(1).
$$
for all
$w \in W^1_{\tau}(0, 1)$.
Multiplying the equation by $(1-i)$, we obtain
\begin{equation}
\matrix{(1-i)\Int_{0}^{1}\lambda \psi {\overline w}\, dx+(1-i)\rho (\lambda+\wp)^{\tilde{\alpha}-1}\psi(1){\overline w}(1)
+(1+i)\Int_{0}^{1} \tau(x) \psi_x {\overline w}_x\, dx=(1-i)\Int_{0}^{1} f_1{\overline w}\, dx\hfill &\cr
-\zeta(1-i)\Int_{-\infty}^{+\infty}\Frac{\eta(\xi)f_2(\xi)}{\lambda +\xi^{2}+\wp}\, d\xi \overline{w}(1)
\ \forall w\in W_{\tau}^1(0, 1).&\cr}
\label{eq66}
\end{equation}
Problem (\ref{eq66}) is of the form
\begin{equation}
{\cal L}(\psi, w)={\cal M}(w),
\label{e26}
\end{equation}
where ${\cal L}:[W_{\tau}^1(0, 1)\times W_{\tau}^1(0, 1)]\longrightarrow \C$
is the bilinear form defined by
$$
{\cal L}(\psi,w)=(1-i)\Int_{0}^{1}\lambda \psi{\overline w}\, dx+(1-i)\rho (\lambda+\wp)^{\tilde{\alpha}-1} \psi(1){\overline w}(1) +(1+i)\Int_{0}^{1} \tau(x) \psi_x {\overline w}_x\, dx,
$$
and ${\cal M}:W_{\tau}^1(0, 1)\longrightarrow \C$
is the antilinear form given by
$$
{\cal M}(w)=(1-i)\Int_{0}^{1} f_1{\overline w}\, dx-\zeta(1-i)\Int_{-\infty}^{+\infty}\Frac{\eta(\xi)f_2(\xi)}{\lambda +\xi^{2}+\wp}\, d\xi \overline{w}(1).
$$
It is easy to verify that ${\cal L}$ is continuous and coercive, and ${\cal M}$ is continuous.
Consequently, by the Lax-Milgram Lemma, system (\ref{e26}) has a unique solution
$\psi\in W_{\tau}^1(0, 1)$.
By the regularity theory for the linear elliptic equations, it follows that $\psi\in H_{\tau}^{2}(0,1)$.
Therefore, the operator $\lambda I- {\cal A}$ is surjective for any $\lambda > 0$. Consequently, using Hill-Yosida
theorem, the result of Theorem \ref{1th} follows.

\section{Asymptotic stability}
In order to state and prove our stability results, we need the following theorem.
\begin{theorem}[\cite{AB}]
Let ${\cal A}$ be the generator of a uniformly bounded $C_0$.
semigroup $\{S(t)\}_{t\geq 0}$ on a Hilbert space ${\cal H}$. If:
\begin{itemize}
\item[(i)] ${\cal A}$ does not have eigenvalues on $i\R$.
\item[(ii)] The intersection of the spectrum $\sigma({\cal A})$ with $i\R$ is at most a countable set,
\end{itemize}
then the semigroup $\{S(t)\}_{t\geq 0}$ is asymptotically stable, i.e, $\|S(t)z\|_{{\cal H}}\rightarrow 0$
as $t\rightarrow \infty$ for any $z\in {\cal H}$.
\label{thm14}
\end{theorem}

\subsection{Strong stability of the system}
Since ${\cal A}$ generates a contraction semigroup and using the Arendt-Batty
Theorem (Theorem \ref{thm14}), system (\ref{10'}) is strongly stable if ${\cal A}$
does not have pure imaginary eigenvalues and $\sigma({\cal A})\cap i\R$ is a countable set.

Our main result is the following theorem:
\begin{theorem}
The $C_0$-semigroup $e^{t{\cal A}}$ is strongly stable in ${\cal H}$; i.e, for all
$U_0\in {\cal H}$, the solution of (\ref{10'}) satisfies
$$
\Lim_{t\rightarrow \infty}\|e^{t{\cal A}}U_0\|_{\cal H}=0.
$$
\label{thmba}
\end{theorem}
For the proof of Theorem \ref{thmba}, we need the following two lemmas.
\begin{lemma}
${\cal A}$ does not have eigenvalues on $i\R$.
\end{lemma}
{\bf Proof}\\
We will argue by contraction. Let us suppose
that there $\lambda\in \R$.\\
$\bullet${\bf Case 1}: $\lambda>0$ and $U\not=0$, such that ${\cal A}U=i\lambda U$. Then, we get
\begin{equation}
\left\{\matrix{i\lambda \psi-\imath (\tau(x) \psi_x)_x=0,\hfill &\cr
i\lambda\vartheta+(\xi^2+\wp)\vartheta-\psi(1)\eta(\xi)=\hfill &\cr}\right.
\label{06z6}
\end{equation}
Then, from (\ref{e17}) and $(\ref{06z6})_2$
we have
\begin{equation}
\vartheta\equiv 0,\ \psi(1)= 0.
\label{e28}
\end{equation}
From the boundary conditions, we have
\begin{equation}
\psi_x(1)=0.
\label{e29}
\end{equation}
We get
\begin{equation}
\left\{\matrix{\lambda \psi-(\tau(x) \psi_x)_x=0,  \hfill & \cr
\psi(1)=\psi_{x}(1)=0, \hfill & \cr
}\right.
\label{4111m}
\end{equation}
Multiplying the equation $(\ref{4111m})_1$ by $\overline{\psi}$
integrating over $(0,1)$ and integrating par parts we get
\begin{equation}\label{mama}
\lambda\|\psi\|^2_{L^2(0,1)}+\int_{0}^{1}\tau(x)|\psi_x|^2dx=0.
\end{equation}
We deduce that
$$
\psi=0.
$$
$\bullet${\bf Case 2}: $\lambda<0$ and $U\not=0$, such that ${\cal A}U=i\lambda U$. Then,
Multiplying $(\ref{4111m})_1$ by $-2x\overline{\psi}_x$, integrating over $(0,1)$, integrating par parts and using the boundary conditions, we obtain
\begin{equation}
\lambda\|\psi\|^2_{L^2(0,1)}-\int_{0}^{1}\tau(x)|\psi_x|^2dx
+\int_{0}^{1}x\tau'(x)|\psi_x|^2dx=0.
\label{eq77}
\end{equation}
Multiplying equations (\ref{mama}) by $-m_{\tau}/2$, and tacking the sum of this equation and (\ref{eq77}), we
get
\begin{equation}
\Frac{2-m_{\tau}}{2}\lambda\Int_{0}^{1}|\psi|^2\, dx-\Int_{0}^{1}\left(\tau(x)-x\tau'(x)+\Frac{m_{\tau}}{2}\tau(x)\right)|\psi_x|^2\, dx=0.
\label{ee92}
\end{equation}
We deduce that
$$
\psi=0.
$$
$\bullet${\bf Case 3}: $\lambda=0$ and $U\not=0$, such that ${\cal A}U=0$.
We get
\begin{equation}
\left\{\matrix{(\tau(x) \psi_x)_x=0,  \hfill & \cr
\psi(1)=\psi_{x}(1)=0, \hfill & \cr
}\right.
\label{4111mm}
\end{equation}
Then $\psi_x= 0$ on $(0, 1)$ and hence $\psi\equiv 0$.
\begin{lemma}
We have
$$
i\R\subset \rho({\cal A}).
$$
\label{leres}
\end{lemma}
{\bf Proof.}

\noindent
$\bullet${\bf Case 1}: $\lambda\not=0$.\\
We will prove that the operator $i\lambda I-{\cal A}$ is surjective for $\lambda\not=0$. For this purpose, let
$F=(f_1, f_2)^{T}\in {\cal H}$, we seek $U=(\psi, \vartheta)^{T}\in D({\cal A})$ solution of the following equation
\begin{equation}
(i\lambda I-{\cal A})U=F.
\label{e35eee}
\end{equation}
Equivalently, we have
\begin{equation}
\left\{\matrix{i\lambda \psi-i(\tau(x) \psi_x)_x=f_1, \hfill & \cr
i\lambda \vartheta+(\xi^{2}+\wp)\vartheta-\psi(1)\eta(\xi)=f_2.\hfill & \cr
}\right.
\label{e35e}
\end{equation}
Solving system (\ref{e35e}) is equivalent to finding $\psi\in H_{\tau}^{2}\cap W_{\tau}^{1}(0,1)$
such that
\begin{equation}
\Int_{0}^{1}(\lambda \psi {\overline w}- (\tau(x) \psi_x)_x {\overline w})\, dx=-i\Int_{0}^{1}f_1{\overline w}\, dx
\label{mk1}
\end{equation}
for all $w\in W_{\tau}^{1}(0,1)$. Then, we get
\begin{equation}
\left\{\matrix{
\Int_{0}^{1}(\lambda \psi {\overline w}+ (\tau(x) \psi_x)  {\overline w_x})\, dx
-i \rho(i\lambda+\wp)^{\tilde{\alpha}-1} \psi(1)\ {\overline w(1)}\, \hfill &\cr
=-i\Int_{0}^1 f_1{\overline w}\, dx+i\zeta \Int_{-\infty}^{+\infty} \Frac{\eta(\xi)}{\xi^{2}+\wp+i\lambda}f_2(\xi){\overline w}  \, d\xi
\,. &\cr}\right.
\label{mk2}
\end{equation}
Pay attention that
$$
-i \rho(i\lambda+\wp)^{\tilde{\alpha}-1} =-i\zeta\Int_{-\infty}^{+\infty} \Frac{\eta(\xi)^2}{(\xi^2+\wp)^2+\lambda^2}\, d\xi
-\zeta\lambda\Int_{-\infty}^{+\infty} \Frac{\eta(\xi)^2}{(\xi^2+\wp)^2+\lambda^2}\, d\xi.
$$
If $\lambda< 0$, we can rewrite (\ref{mk2}) as
\begin{equation}
-({\cal L}_{\lambda}\psi, w)_{W_{\tau}^1}+(\psi, w)_{W_{\tau}^1}={\cal M}(w)
\label{e2xx5}
\end{equation}
with the inner product defined by
$$
\matrix{(\psi, w)_{W_{\tau}^1}=\Int_{0}^{1}  \tau(x)\psi_x {\overline w_x}\, dx
-\zeta\lambda\Int_{-\infty}^{+\infty} \Frac{\eta(\xi)^2}{(\xi^2+\wp)^2+\lambda^2}\, d\xi \psi(1)\ {\overline w(1)}\hfill & \cr
\qquad -i\zeta\Int_{-\infty}^{+\infty} \Frac{\eta(\xi)^2}{(\xi^2+\wp)^2+\lambda^2}\, d\xi \psi(1)\ {\overline w(1)}&\cr}
$$
and
$$
({\cal L}_{\lambda}\psi, w)_{W_{\tau}^1}=-\Int_{0}^1\lambda \psi {\overline w}\, dx.
$$
Using the compactness embedding from $L^2(0,1)$ into $W_{\tau}^{-1}(0,1)$ and from $W_{\tau}^1(0,1)$ into $L^2(0,1)$
we deduce that the operator ${\cal L}_\lambda$ is compact from $L^2(0,1)$ into $L^2(0,1)$. Consequently, by Fredholm
alternative, proving the existence of $\psi$ solution of (\ref{e35e}) reduces to proving that $1$ is not an eigenvalue
of ${\cal L}_{\lambda}$. Indeed if $1$ is an eigenvalue, then there exists $\psi\not=0$, such that
\begin{equation}
({\cal L}_{\lambda}\psi, w)_{W_{\tau}^1}=(\psi, w)_{W_{\tau}^1}\quad \forall w\in W_{\tau}^1.
\label{eaz}
\end{equation}
In particular for $w=\psi$, it follows that
$$
-\lambda\|\psi(x)\|_{L^2(0,1)}^2+i\rho(i\lambda+\wp)^{\tilde{\alpha}-1} |\psi(1)|^2
=\|\sqrt{\tau(x)}\psi_x(x)\|_{L^2(0,1)}^{2}.
$$
Hence, we have
\begin{equation}
\psi(1)=0.
\label{e22xx1}
\end{equation}
From (\ref{eaz}), we obtain
\begin{equation}
(\tau\psi_x)(1)=0
\label{e23xx1}
\end{equation}
Then
\begin{equation}
\left\{\matrix{\lambda \psi-(\tau(x) \psi_x)_x=0 \hbox { on } (0,1), \hfill & \cr
\left\{\matrix{\psi(0)=0 \hfill & \hbox{ if } m_{\tau}\in [0, 1)\hfill &\cr
(\tau(x) \psi_x)(0)=0\hfill & \hbox{ if } m_{\tau}\in [1, 2)\hfill &\cr}\right.\hfill & \cr
\psi(1)=0,\, \psi_x(1)=0.  \hfill & \cr}\right.
\label{e2xx1}
\end{equation}
We deduce that $U=0$. \\

If $\lambda> 0$, we can rewrite (\ref{mk2}) as
\begin{equation}
-({\cal L}_{\lambda}\psi, w)_{W_{\tau}^1}+(\psi, w)_{W_{\tau}^1}={\cal M}(w),
\label{e2xx59}
\end{equation}
with the inner product defined by
$$
(\psi, w)_{W_{\tau}^1}=\Int_{0}^{1}  \tau(x)\psi_x {\overline w_x}\, dx+\Int_{0}^1\lambda \psi {\overline w}\, dx
-i\zeta\Int_{-\infty}^{+\infty} \Frac{\eta(\xi)^2}{(\xi^2+\wp)^2+\lambda^2}\, d\xi \psi(1)\ {\overline w(1)}.
$$
and
$$
({\cal L}_{\lambda}\psi, w)_{W_{\tau}^1}=\zeta\lambda\Int_{-\infty}^{+\infty} \Frac{\eta(\xi)^2}{(\xi^2+\wp)^2+\lambda^2}\, d\xi \psi(1)\ {\overline w(1)}.
$$
In a similar way, we prove that $U=0$.

\noindent
$\bullet${\bf Case 2}: $\lambda =0$ and $\wp\not=0$.\\ The system (\ref{e35e}) is reduced to the following system
\begin{equation}
\left\{\matrix{-(\tau(x) \psi_x)_x=-if_1, \hfill & \cr
(\xi^{2}+\wp)\vartheta-\psi(1)\eta(\xi)=f_2.\hfill & \cr
}\right.
\label{e40e}
\end{equation}
Solving system (\ref{e40e}) is equivalent to finding $\psi\in W_{\tau}^{1}(0,1)$
such that
\begin{equation}
-\Int_{0}^{1}(\tau(x) \psi_x)_x {\overline w}\, dx=-i\Int_{0}^{1}f_1{\overline w}\, dx
\label{mk10}
\end{equation}
for all $w\in W_{\tau}^{1}(0,1)$. Then, we get
\begin{equation}
\Int_{0}^{1}\tau(x) \psi_x {\overline w_x}\, dx-i\rho \wp^{\alpha-1}\psi(1) {\overline w(1)}
=-i\Int_{0}^{1}f_1{\overline w}\, dx +i\zeta \Int_{-\infty}^{\infty}\Frac{\eta(\xi)f_2 (\xi)}{\xi^2+\wp}\, d\xi
{\overline w(1)}.
\label{mk11za}
\end{equation}
Multiplying this equation by $(1+i)$, we obtain
\begin{equation}
\matrix{
(1+i)\Int_{0}^{1}\tau(x) \psi_x {\overline w_x}\, dx-(1+i)i\rho \wp^{\alpha-1}\psi(1) {\overline w(1)}
=-i(1+i)\Int_{0}^{1}f_1{\overline w}\, dx\hfill  &\cr
\hspace{4cm} +i(1+i)\zeta \Int_{-\infty}^{\infty}\Frac{\eta(\xi)f_2 (\xi)}{\xi^2+\wp}\, d\xi {\overline w(1)}.\hfill  &\cr}
\label{mk11}
\end{equation}
Consequently, (\ref{mk11}) is equivalent to the problem
\begin{equation}
{\cal L}(\psi, w)={\cal M}(w),
\label{mk12}
\end{equation}
where the bilinear form ${\cal L}: W_{\tau}^{1}(0,1)\times W_{\tau}^{1}(0,1)\rightarrow \C$
and the linear form ${\cal M}:W_{\tau}^{1}(0,1)\rightarrow \C$ are defined by
\begin{equation}
{\cal L}(\psi, w)=(1+i)\Int_{0}^{1}\tau(x) \psi_x {\overline w_x}\, dx-(1+i)i\rho \wp^{\alpha-1}\psi(1) {\overline w(1)}
\label{mk13}
\end{equation}
and
$$
{\cal M}(w)=-i(1+i)\Int_{0}^{1}f_1{\overline w}\, dx
+i(1+i)\zeta \Int_{-\infty}^{\infty}\Frac{\eta(\xi)f_2 (\xi)}{\xi^2+\wp}\, d\xi {\overline w(1)}.
$$
It is easy to verify that ${\cal L}$ is continuous and coercive, and ${\cal M}$ is continuous. So by applying
the Lax-Milgram theorem, we deduce that for all $w\in W_{\tau}^{1}(0,1)$
problem (\ref{mk12}) admits a unique solution $\psi\in W_{\tau}^{1}(0, 1)$.
Applying the classical elliptic regularity, it follows from (\ref{mk11}) that $\psi\in H_{\tau}^{2}(0,1)$.
Therefore, the operator ${\cal A}$ is surjective.

\hfill$\diamondsuit$\\

\subsection{Stability}
In this section, we will study the exponential decay of solution of the system $(\ref{10'})$. For this purpose we will
use the following theorem.
Next, in order to prove an polynomial decay rate we will use the following theorem.

\begin{theorem}\label{5.2}{\bf\cite{BT}}
Let ${\cal A}$ be the generator of a strongly continuous
semigroup of contractions $(e^{t{\cal A}})_{t\geq 0}$ on a Hilbert space ${\cal X}$. If
\begin{equation}
\imath \R\subset\rho({\cal A}),
\label{aazzss}
\end{equation}
where $\rho({\cal A})$ is the resolvent set of ${\cal A}$, then for fixed $l > 0$, the following conditions
are equivalent:
$$
\overline{\lim_{|\beta|\rightarrow\infty}}\frac{1}{\beta^l}
\|(\imath\beta I-{{\cal A}})^{-1}\|_{{\cal L}({\cal X})}<\infty
$$
for some $l$, then there exist $c$ such that
$$
\|e^{{\cal A}t}V_0\|^2\leq \frac{c}{t^{\frac{2}{l}}}\|V_0\|^2_{D({\cal A})}.
$$
\end{theorem}

Our main result is the following.
\begin{theorem}\label{4.2}
If $\wp\not=0$, then The semigroup $({\cal S}_{\cal A}(t))_{t\geq 0}$ is polynomially stable and
$$
E(t)=\|S_{\cal A}(t)V_0\|_{\cal H}^2\leq
\Frac{c}{ t^{\frac{2}{1-\tilde{\alpha}}}}\|V_0\|_{D({\cal A})}^2.
$$
Moreover, the rate of energy decay $t^{\frac{2}{1-\tilde{\alpha}}}$ is sharp for general initial data in $D({\cal A})$.
\end{theorem}
{\bf Proof.}
We will need to study the resolvent equation $(i\lambda-{\cal A})U=F$, for $\lambda\in \R$, namely
\begin{equation}
\left\{\matrix{i\lambda \psi-\imath (\tau(x) \psi_x)_x=f_1,  \hfill & \cr
i\lambda \vartheta+(\xi^{2}+\wp)\vartheta-\psi(1)\eta(\xi)=f_2.\hfill & \cr}\right.
\label{e188}
\end{equation}
We aim to show that there exist $C$ such that
$$
\|(\imath\lambda I-{\cal A})^{-1}\|_{\mathcal{L(H)}}<C|\lambda|^{1-\tilde{\alpha}}.
$$
 Taking inner product in ${\cal H}$ with $u$
and using (\ref{e17}) we get
\begin{equation}
|Re\langle{\cal A}U, U\rangle|\leq \|U\|_{\cal H}\|F\|_{\cal H}.
\label{e360e1}
\end{equation}
This implies that
\begin{equation}
\zeta \Int_{-\infty}^{+\infty}(\xi^{2}+\wp)|\vartheta(\xi)|^2d\xi\leq \|U\|_{\cal H}\|F\|_{\cal H}.
\label{e37}
\end{equation}
By multiplying $(\ref{e188})_2$ by $(i\lambda+\xi^{2}+\wp)^{-2}|\xi|$, we get
\begin{equation}
(i\lambda+\xi^{2}+\wp)^{-2} \psi(1)\eta(\xi)|\xi|=(i\lambda+\xi^{2}+\wp)^{-1}|\xi|\vartheta-(i\lambda+\xi^{2}+\wp)^{-2}|\xi|f_2(\xi).
\label{e38kk}
\end{equation}
Hence, by taking absolute values of both sides of (\ref{e38kk}), integrating over the interval $]-\infty, +\infty[$ with
respect to the variable $\xi$ and applying Cauchy-Schwartz inequality, we obtain
\begin{equation}
{\cal R}|\psi(1)|\leq  \sqrt{2}{\cal P} \left(\Int_{-\infty}^{+\infty}\xi^{2}|\vartheta|^{2}\, d\xi\right)^{\frac{1}{2}}
+ 2 {\cal Q}\left(\Int_{-\infty}^{+\infty}|f_2(\xi)|^{2}\, d\xi\right)^{\frac{1}{2}},
\label{e39kk}
\end{equation}
where
$$
{\cal R}=\left|\Int_{-\infty}^{+\infty}(i\lambda+\xi^{2}+\wp)^{-2}|\xi|\eta(\xi)\, d\xi\right|
=\frac{|1-2\tilde{\alpha}|}{4}\Frac{\pi}{|\sin\frac{(2\tilde{\alpha}+3)}{4}\pi|}|i\lambda+\wp|^{\frac{(2\tilde{\alpha}-5)}{4}},
$$
$$
{\cal P}=\left(\Int_{-\infty}^{+\infty}(|\lambda|+\xi^{2}+\wp)^{-2}\, d\xi\right)^{\frac{1}{2}}=(\frac{\pi}{2})^{1/2}||\lambda|+\wp|^{-\frac{3}{4}},
$$
$$
{\cal Q}=\left(\Int_{-\infty}^{+\infty}(|\lambda|+\xi^{2}+\wp)^{-4}|\xi|^{2}\, d\xi\right)^{\frac{1}{2}}
=\left(\Frac{\pi}{16}||\lambda|+\wp|^{-\frac{5}{2}}\right)^{1/2}.
$$
We deduce that
\begin{equation}
|\psi(1)|^2\leq c|\lambda|^{1-\tilde{\alpha}}\|U\|_{\cal H}\|F\|_{\cal H}+c\|F\|_{\cal H}^{2}.
\label{e41nkk}
\end{equation}
We will use the classical multiplier method. Let us multiply the equation $(\ref{e188})$ by $\overline{\psi}$
integrating over $(0,1)$ and integrating par parts we get
\begin{equation}\label{ma}
\lambda\|\psi\|^2_{L^2(0,1)}-\left[\tau(x)\psi_x\overline{\psi}\right]_0^1+\int_{0}^{1}\tau(x)|\psi_x|^2dx=-\imath\int_{0}^{1}f_1\overline{\psi}dx.
\end{equation}
If $\lambda> 0$ large, we have from (\ref{ma}) and  (\ref{e41nkk}) that
$$
\matrix{\lambda\|\psi\|^2_{L^2(0,1)}&\leq & \|f_1\|_{L^2(0, 1)}\|\psi\|_{L^2(0, 1)}+|(\tau(x)\psi_x)(1)||\overline{\psi}(1)|\hfill \cr
&\leq & \|f_1\|_{L^2(0, 1)}\|\psi\|_{L^2(0, 1)}+\zeta|\Int_{\R}\eta(\xi)\vartheta(\xi)\, d\xi||\overline{\psi}(1)|\hfill \cr
&\leq & \|F\|_{\cal H}\|U\|_{\cal H}+c(\|F\|_{\cal H}\|U\|_{\cal H})^{1/2}(|\lambda|^{\frac{1-\tilde{\alpha}}{2}}(\|U\|_{\cal H}\|F\|_{\cal H})^{1/2}+c\|F\|_{\cal H}).\hfill \cr
&\leq & c\|F\|_{\cal H}\|U\|_{\cal H}+c|\lambda|^{\frac{1-\tilde{\alpha}}{2}}\|F\|_{\cal H}\|U\|_{\cal H}+c\|F\|_{\cal H}^2.\hfill \cr
}
$$
Then, we deduce that
$$
\|U\|_{\cal H}^2\leq \Frac{C}{\lambda} \|F\|_{\cal H}^2.
$$
Now, if $\lambda< 0$ large,
we multiply (\ref{e188}) by $-2x\overline{\psi}_x$, integrating over $(0,1)$, and integrating par parts we obtain
\begin{eqnarray*}
-\lambda\left[x|\psi|^2\right ]_0^1&+&\lambda\|\psi\|^2_{L^2(0,1)}+2\left[\tau(x)|\psi_x|^2x\right]_0^1-2\int_{0}^{1}\tau(x)|\psi_x|^2dx-\left[x\tau(x)|\psi_x|^2\right]_0^1\\\\&+&
\int_{0}^{1}\tau(x)|\psi_x|^2dx+\int_{0}^{1}x\tau'(x)|\psi_x|^2dx=2\Re\ i\int_{0}^{1}f_1 x\overline{\psi}_xdx.
\end{eqnarray*}
Multiplying $(\ref{ma})$ by $-\frac{m_{\tau}}{2}$, summing with the last equation and taking the real part, we get
\begin{eqnarray*}
\lambda\left(1-\frac{m_{\tau}}{2}\right)\|\psi\|^2_{L^2(0,1)}&-&\int_{0}^{1}\left(\tau(x)-x\tau'(x)
+\frac{m_{\tau}}{2}\tau(x)\right)|\psi_x|^2dx-\lambda\left[x|\psi^2|\right]_0^1+\left[x\tau(x)|\psi_x|^2\right]_0^1\\\\&+&
\frac{m_{\tau}}{2}\left[\tau(x)\psi_x\overline{\psi}\right]_0^1=2\imath\Re\int_0^1f_1 x\overline{\psi}_xdx
+\frac{m_{\tau}}{2}\imath\int_0^1f_1\overline{\psi}dx.
\end{eqnarray*}
Moreover
$$
\Int_{0}^{1}x^2 |\psi_x|^2\, dx\leq \Int_{0}^{1}x^{m_{\tau}} |\psi_x|^2\, dx\leq \Frac{1}{\tau(1)}\Int_{0}^{1}\tau(x) |\psi_x|^2\, dx
$$
Then
$$
\matrix{\left|\Int_{0}^{1}f_1 x \overline{\psi}_x\, dx\right|&\leq &\|f_1\|_{L^2(0, 1)}\|x \psi_x\|_{L^2(0, 1)}\hfill \cr
&\leq &\Frac{1}{\sqrt{\tau(1)}}\|f_1\|_{L^2(0, 1)}\|\sqrt{\tau(x)} \psi_x\|_{L^2(0, 1)}\hfill \cr
&\leq &c(\varepsilon)\|f_1\|_{L^2(0, 1)}^2+\varepsilon\|\sqrt{\tau(x)} \psi_x\|_{L^2(0, 1)}^2\hfill \cr
&\leq &c(\varepsilon)\|F\|_{\cal H}^2+\varepsilon\|\sqrt{\tau(x)} \psi_x\|_{L^2(0, 1)}^2\hfill \cr}
$$
Then, we deduce that
\begin{equation}
\matrix{-\lambda\left(1-\frac{m_{\tau}}{2}\right)\|\psi\|^2_{L^2(0,1)}+\Int_{0}^{1}\left(\tau(x)-x\tau'(x)+\frac{m_{\tau}}{2}\tau(x)\right)|\psi_x|^2dx
=-\lambda |\psi(1)|^2  \hfill &\cr
+\Frac{1}{\tau(1)}|\tau(1)\psi_x(1)|^2+\frac{ m_{\tau}}{2}(\tau(1)\psi_x(1)\overline{\psi}(1))-2\imath\Re\int_0^1f_1 x\overline{\psi}_xdx
-\frac{m_{\tau}}{2}\imath\int_0^1f_1\overline{\psi}dx.\hfill &\cr}
\label{eq44}
\end{equation}

Then, there exist $c_1$, $c_2$ and $c_3$ such that
$$
\matrix{|\lambda|\left(1-\frac{m_{\tau}}{2}\right)\|\psi\|^2_{L^2(0,1)}+\int_{0}^{1}\left(\tau(x)-x\tau'(x)+\frac{m_{\tau}}{2}\tau(x)\right)|\psi_x|^2dx \leq \hfill &\cr
c|\lambda|(|\lambda|^{1-\tilde{\alpha}}\|U\|_{\cal H}\|F\|_{\cal H}+c\|F\|_{\cal H}^{2})
+c(\|F\|_{\cal H}\|U\|_{\cal H})^{1/2}(|\lambda|^{\frac{1-\tilde{\alpha}}{2}}(\|U\|_{\cal H}\|F\|_{\cal H})^{1/2}+c\|F\|_{\cal H})\hfill &\cr
+c\|U\|_{\cal H}\|F\|_{\cal H}+c\|F\|_{\cal H}\|\sqrt{\tau(x)} \psi_x\|_{L^2(0, 1)}\hfill &\cr}
$$
By definition of $m_{\tau}$, we have
$$
\Frac{(2-m_{\tau})}{2}\tau\leq (\tau-x\tau')+\Frac{m_{\tau}}{2} \tau
$$
Then
$$
\matrix{\|\psi\|^2_{L^2(0,1)}\leq c |\lambda|^{1-\tilde{\alpha}}\|U\|_{\cal H}\|F\|_{\cal H}
+c|\lambda|^{-\frac{1+\tilde{\alpha}}{2}}\|U\|_{\cal H}\|F\|_{\cal H}\hfill &\cr
+\Frac{c}{|\lambda|}(\|U\|_{\cal H}\|F\|_{\cal H})^{1/2}\|F\|_{\cal H}+\Frac{c}{|\lambda|}\|U\|_{\cal H}\|F\|_{\cal H}
+\Frac{c}{|\lambda|}\|F\|_{\cal H}^2+c\|F\|_{\cal H}^2\hfill &\cr
\leq c |\lambda|^{1-\tilde{\alpha}}\|U\|_{\cal H}\|F\|_{\cal H}+c'\|F\|_{\cal H}^2.\hfill &\cr}
$$
Since that ($\wp> 0$)
$$
\Int_{-\infty}^{+\infty}|\vartheta(\xi)|^{2}\, d\xi
\leq \Frac{1}{\wp}  \Int_{-\infty}^{+\infty}(\xi^{2}+\wp)|\vartheta(\xi)|^{2}\, d\xi\leq C\|U\|_{\cal H}\|F\|_{\cal H}.
$$
Hence
\begin{equation}
\|U\|_{\cal H}^2\leq c |\lambda|^{1-\tilde{\alpha}}\|U\|_{\cal H}\|F\|_{\cal H}.
\label{ee25}
\end{equation}
Then
$$
\|U\|_{\cal H}\leq c |\lambda|^{1-\tilde{\alpha}}\|F\|_{\cal H}.
$$
Hence, for large $\lambda$
$$
\|(\imath\lambda I-{\cal A})^{-1}F\|_{\mathcal{H}}<C |\lambda|^{1-\tilde{\alpha}}\|F\|_{L^2(0,1)}.
$$
Therefore
$$
\|(\imath\lambda I-{\cal A})^{-1}\|_{\mathcal{L(H)}}<C |\lambda|^{1-\tilde{\alpha}}.
$$
The conclusion then follows by applying the Theorem \ref{5.2}.
\section{Optimality of energy decay}
In this section, we will study the lack of exponential decay of solution of the system $(\ref{10'})$. For this purpose
we will use the following theorem.
\begin{lemma}\label{4.1}{\bf\cite{P}}
Let $S(t)$ be a $C_0$-semigroup of contractions on Hilbert space $\mathcal{X}$ with generator ${\mathcal{A}}$. Then $S(t)$ is exponentially stable if and only if
$$
\rho({\mathcal{A}})\supseteq\{\imath\beta: \beta \in \R\}\equiv\imath\R
$$
and
$$
\overline{\lim_{|\beta|\rightarrow\infty}}
\|(\imath\beta I-{\mathcal{A}})^{-1}\|_{\mathcal{L(X)}}<\infty.
$$
\end{lemma}

Our main result is the following.
\begin{theorem}\label{4.28}
The semigroup generated by the operator ${\mathcal{A}}$ is not exponentially stable.
\end{theorem}
{\bf Proof.} We consider the case $\tau(x)=x^{\alpha}, 0<\alpha<1$.
We aim to show that an infinite number of eigenvalues of ${\mathcal{A}}$ approach the imaginary axis
which prevents the system $(\ref{10'})$ from being exponentially stable. Let $\lambda$ be an eigenvalue of
${\mathcal{A}}$ with associated eigenvector $U$. Then the equation ${\mathcal{A}}U=\lambda U$ is equivalent to
$$
\imath\lambda \psi+(x^\alpha \psi_x)_x=0
$$
together with the conditions
$$
\quad \left \{
\begin{array}{ll}
\psi_x(1)=\imath\rho (\lambda+\wp)^{\tilde{\alpha}-1}\psi(1),\\
\psi(0)=0,
\end{array}
\right.
$$
so we get the following system
\begin{equation}\label{40}
\quad \left \{
\begin{array}{ll}
\gamma^2\psi-(x^\alpha \psi_x)_x=0,\\
\psi_x(1)=\imath\rho (\lambda+\wp)^{\tilde{\alpha}-1}\psi(1),\\
\psi(0)=0,
\end{array}
\right.
\end{equation}
with $\gamma^2=-\imath\lambda$.\\

Suppose that $v$ is a solution of $(\ref{40})_1$, then the function $\Psi$ defined by
$$ 
v(x)=x^{\frac{1-\alpha}{2}}\Psi\left(\frac{2}{2-\alpha} \imath\gamma
x^{\frac{2-\alpha}{2}} \right)
$$ 
is a solution of the following equation
\begin{equation}\label{41}
y^2\Psi''(y)+y\Psi'(y)+\left(y^2-\left( \frac{\alpha-1}{2-\alpha}\right) ^2 \right) \Psi(y)=0.
\end{equation}
We have
\begin{equation}\label{42}
\psi(x)=c_+\widetilde{\theta}_+(x)+c_-\widetilde{\theta}_-(x)
\end{equation}
where
$$
\widetilde{\theta}_+(x)=x^{\frac{1-\alpha}{2}}J_{\nu_{\alpha}}\left(\frac{2}{2-\alpha}\imath\gamma x^{\frac{2-\alpha}{2}}\right)
\quad \hbox{ and  } \ \ \
\widetilde{\theta}_-(x)=x^{\frac{1-\alpha}{2}}J_{-\nu_{\alpha}}\left(\frac{2}{2-\alpha}\imath\gamma x^{\frac{2-\alpha}{2}}\right).
$$
$\psi(0)=0$, then $c_{-}{\tilde d}^ {-}=0$, then $c_{-}=0$. So
$$
\psi(x)=c_+\widetilde{\theta}_+(x).
$$
Therefore the boundary conditions can be written as the following system
$$
c_+\widetilde{\theta}'_+(1)=\imath\rho (\lambda+\wp)^{\tilde{\alpha}-1}c_+\widetilde{\theta}_+(1).
$$
Then, a non-trivial solution $u$ exists if and only if
$$
f(\gamma)=(1-\alpha)
J_{\nu_\alpha}\left(\frac{2\gamma}{2-\alpha}\imath\right)-
\imath \gamma
J_{1+\nu_\alpha}\left(\frac{2\gamma}{2-\alpha}\imath\right)
-i\rho(\lambda+\wp)^{\tilde{\alpha}-1}
J_{\nu_\alpha}\left(\frac{2\gamma}{2-\alpha}\imath\right)=0.
$$
Our purpose is to prove, thanks to Rouch\'{e}'s Theorem, that there is a subsequence of eigenvalues for which their
real part tends to $0$.

In the sequel, since ${\mathcal{A}}$ is dissipative, we study the asymptotic behavior of the large eigenvalues $\lambda$
of ${\mathcal{A}}$ in the strip $-\alpha_0\leq\Re(\lambda)\leq0$, for some $\alpha_0>0$ large enough and for such
$\lambda$, we remark that $\theta_+$ and $\theta_-$ remain bounded.
\begin{lemma}\label{l4.1}
There exists $N \in \N$ such that
$$
\{\lambda_k\}_{{k\in\Z}^*, \ |k|\geq N} \subset \sigma({\mathcal{A}}),
$$
where

\noindent
$\bullet$ If $1\geq\tilde{\alpha}> \frac{1}{2}$, then
\begin{eqnarray*}
\lambda_k&=&-i \left[C_0^2(k\pi)^2+C_1^2\pi^2+2 C_0 C_1 k\pi^2+2{C_0 C_2}\right]
-2i^{1-\tilde{\alpha}}\Frac{C_0 C_3}{(k\pi)^{2-2\tilde{\alpha}}}+o\left(\frac{1}{k^{2-2\tilde{\alpha}}}\right),
\end{eqnarray*}
where
$$
C_0=-\frac{2-\alpha}{2}, C_1=-\frac{2-\alpha}{2}\left(\frac{\nu_{\alpha}}{2}+\frac{5}{4}\right),
$$
$$
C_2=\Frac{(2-\alpha)}{4}\left((\frac{1}{2}+\nu_{\alpha})(\frac{3}{2}+\nu_{\alpha})-\frac{4(1-\alpha)}{(2-\alpha)}\right)
$$
and
$$
C_3=-\rho \left(\frac{2}{2-\alpha}\right)^{2-2\tilde{\alpha}}.
$$
\noindent
$\bullet$ If $\tilde{\alpha}< \frac{1}{2}$, then
\begin{eqnarray*}
\lambda_k&=&\imath\gamma_k^2\\ \\
&=&-i \left[C_0^2(k\pi)^2+C_1^2\pi^2+2 C_0 C_1 k\pi^2+2{C_0 C_2}+2\Frac {{C_0 C_3}}{k}
+2\Frac {{C_1 C_2}}{k}\right]\\ \\
&& -2i^{1-\tilde{\alpha}}\Frac{C_0 C_4}{(k\pi)^{2-2\tilde{\alpha}}}+o\left(\frac{1}{k^{2-2\tilde{\alpha}}}\right).
\end{eqnarray*}
where
$$
C_0=-\frac{2-\alpha}{2}, C_1=-\frac{2-\alpha}{2}\left(\frac{\nu_{\alpha}}{2}+\frac{5}{4}\right),
$$
$$
C_2=\Frac{(2-\alpha)}{4}\left((\frac{1}{2}+\nu_{\alpha})(\frac{3}{2}+\nu_{\alpha})-\frac{4(1-\alpha)}{(2-\alpha)}\right)
$$
$$
C_3=-\frac{2-\alpha}{4}\left(-\Frac{(\frac{1}{2}+\nu_{\alpha})(\frac{3}{2}+\nu_{\alpha})}{4}(2-\alpha)+(1-\alpha)\right)
\Frac{2C_1}{C_0^2}.
$$
and
$$
C_4=-\rho \left(\frac{2}{2-\alpha}\right)^{2-2\tilde{\alpha}}.
$$
and
$$
\lambda_k=\overline{\lambda_{-k}}, \hbox{ if } k \leq -N.
$$
\end{lemma}
{\bf Proof.}
We will use the following classical development (see {\bf\cite{Le}}), for all $\delta>0$ and when $|\arg{z}|<\pi-\delta$:
\begin{equation}
\matrix{J_{\nu}(z)=\left(\Frac{2}{\pi z}\right)^{1/2}\left[\cos\left(z-\nu\frac{\pi}{2}-\frac{\pi}{4}\right)
-\Frac{(\nu-\Frac{1}{2})(\nu+\Frac{1}{2})}{2}\Frac{\sin\left(z-\nu\frac{\pi}{2}-\frac{\pi}{4}\right)}{z}\right.\hfill & \cr
\left.-\Frac{(\nu-\Frac{1}{2})(\nu+\Frac{1}{2})(\nu-\Frac{3}{2})(\nu+\Frac{3}{2})}{8}\Frac{\cos\left(z-\nu\frac{\pi}{2}-\frac{\pi}{4}\right)}{z^2}
+O\left(\frac{1}{|z|^3}\right)\right]. & \cr}
\label{84}
\end{equation}
We get
$$
f(\gamma)=-\imath\gamma
\left( \frac{2}{\pi z}\right)^\frac{1}{2}
\frac{e^{-\imath(z-(\nu_{\alpha}+1)\frac{\pi}{2}-\frac{\pi}{4})}}{2}\widetilde{f}(\gamma),
$$
where $$z=\frac{2\gamma}{2-\alpha}\imath$$
and
\begin{eqnarray*}
\widetilde{f}(\gamma)&=&1+e^{2\imath(z-(\nu_{\alpha}+1)\frac{\pi}{2}-\frac{\pi}{4})}-
\Frac{(\frac{1}{2}+\nu_{\alpha})(\frac{3}{2}+\nu_{\alpha})}{2i}\left( \frac{2}{2-\alpha}\imath\right)^{-1}
\Frac{e^{2\imath(z-(\nu_{\alpha}+1)\frac{\pi}{2}-\frac{\pi}{4})}-1}{\gamma}\\ \\
&-&\frac{1-\alpha}{i\gamma}(e^{\imath(2z-\nu_{\alpha}\pi-\pi)}+e^{-i\frac{\pi}{2}})
+\frac{\rho}{\gamma}\lambda^{\tilde{\alpha}-1}(e^{\imath(2z-\nu_{\alpha}\pi-\pi)}+e^{-i\frac{\pi}{2}})\\ \\
&-&\Frac{(\nu_{\alpha}+\Frac{1}{2})(\nu_{\alpha}+\Frac{3}{2})(\nu_{\alpha}-\Frac{1}{2})(\nu_{\alpha}+\Frac{5}{2})}{8}\Frac{1+e^{2\imath(z-(\nu_{\alpha}+1)\frac{\pi}{2}-\frac{\pi}{4})}}{z^2}\\ \\
&-&\frac{1-\alpha}{i\gamma}\Frac{(\frac{1}{2}-\nu_{\alpha})(\frac{1}{2}+\nu_{\alpha})}{2}\Frac{(e^{\imath(2z-\nu_{\alpha}\pi-\pi)}-e^{-i\frac{\pi}{2}})}{i z}\\ \\
&=&f_0(\gamma)+\frac{f_1(\gamma)}{\gamma}+\frac{f_2(\gamma)}{\gamma^{3-2\tilde{\alpha}}}+\frac{f_3(\gamma)}{\gamma^{2}}
+O\left(\frac{1}{\gamma^3} \right)+O\left(\frac{1}{\gamma^{4-2\tilde{\alpha}}} \right)
\end{eqnarray*}
with
$$
\matrix{f_0(\gamma)&=&1+e^{2\imath(z-(\nu_{\alpha}+1)\frac{\pi}{2}-\frac{\pi}{4})}, \hfill \cr
f_1(\gamma)&=&-
\Frac{(\frac{1}{2}+\nu_{\alpha})(\frac{3}{2}+\nu_{\alpha})}{2i}\left( \frac{2}{2-\alpha}\imath\right)^{-1}
(e^{2\imath(z-(\nu_{\alpha}+1)\frac{\pi}{2}-\frac{\pi}{4})}-1)
-\frac{1-\alpha}{i}(e^{\imath(2z-\nu_{\alpha}\pi-\pi)}+e^{-i\frac{\pi}{2}}),\hfill \cr
f_2(\gamma)&=&\rho i^{\tilde{\alpha}-1}(e^{\imath(2z-\nu_{\alpha}\pi-\pi)}+e^{-i\frac{\pi}{2}}).\hfill \cr
f_3(\gamma)&=&-\Frac{(\nu_{\alpha}+\Frac{1}{2})(\nu_{\alpha}+\Frac{3}{2})(\nu_{\alpha}-\Frac{1}{2})(\nu_{\alpha}+\Frac{5}{2})}{8}\left( \frac{2}{2-\alpha}\imath\right)^{-2}(1+e^{2\imath(z-(\nu_{\alpha}+1)\frac{\pi}{2}-\frac{\pi}{4})})\hfill  \cr
&&\ \ \ +\frac{1-\alpha}{\gamma}\Frac{(\frac{1}{2}-\nu_{\alpha})(\frac{1}{2}+\nu_{\alpha})}{2}\left( \frac{2}{2-\alpha}\imath\right)^{-1}(e^{\imath(2z-\nu_{\alpha}\pi-\pi)}-e^{-i\frac{\pi}{2}}),\hfill & \cr}
$$
Note that $f_0, f_1, f_2$ and $f_3$ remain bounded in the strip $-\alpha_0\leq\Re(\lambda)\leq0$.\\

\noindent
{\bf Case 1: $\tilde{\alpha}>\frac{1}{2}$}.\\

\noindent
We search the roots of $f_0$,
$$
f_0(\gamma)=0 \ \Leftrightarrow \ 1+e^{2\imath(z-(\nu_{\alpha}+1)\frac{\pi}{2}-\frac{\pi}{4})}=0,
$$
so, $f_0$ has the following roots
$$
\gamma_k^0=-\frac{2-\alpha}{2}\imath\left(k+\frac{\nu_{\alpha}}{2}+\frac{5}{4}\right)\pi, \ \ k \in \Z.
$$
Let $B_k(\gamma_k^0,r_k)$ be the ball of centrum $\gamma_k^0$ and
radius $r_k=\frac{1}{\sqrt{|k|}}$, then if $\gamma\in \partial B_k$, we have $\gamma=\gamma_k^0+r_ke^{\imath\theta}, \ \theta\in [0,2\pi] $, then we have
$$
f_0(\gamma)=\frac{4}{2-\alpha}r_ke^{\imath\theta}+o(r_k^2).
$$
Hence, there exists a positive constant c such that, for all $\gamma\in \partial B_k$
$$
|f_0(\gamma)|\geq cr_k=\frac{c}{\sqrt{k}}.
$$
From the expression of $\widetilde{f}$, we conclude that
$$
|\widetilde{f}(\gamma)-f_0(\gamma)|=O\left(\frac{1}{\gamma} \right) =O\left(\frac{1}{k} \right),
$$
then, for $k$ large enough, for all $\gamma\in \partial B_k$
$$
|\widetilde{f}(\gamma)-f_0(\gamma)|<|f_0(\gamma)|.
$$

Using Rouch\'{e}'s Theorem, we deduce that $\widetilde{f}$ and $f_0$ have the same number of zeros in $B_k$.
Consequently, there exists a subsequence of roots of $\widetilde{f}$ that tends to the roots $\gamma_k^0$ of $f_0$,
then there exists $N\in {\N}$ and a subsequence $\{\gamma_k\}_{|k|\geq N}$ of roots of $f(\gamma)$, such that
$\gamma_k=\gamma_k^0+o(1)$ that tends to the roots
$-\frac{2-\alpha}{2}\imath\left(k+\frac{\nu_{\alpha}}{2}+\frac{5}{4}\right)\pi$ of $f_0$. 
\\

Now, we can write
\begin{equation}\label{89}
\gamma_k=-\frac{2-\alpha}{2}\imath\left(k+\frac{\nu_{\alpha}}{2}+\frac{5}{4}\right)\pi+\varepsilon_k,
\end{equation}
then
\begin{eqnarray*}
e^{2\imath(z-(\nu_{\alpha}+1)\frac{\pi}{2}-\frac{\pi}{4})}&=&-e^{-\frac{4}{2-\alpha}\varepsilon_k}\\
&=&-1+\frac{4}{2-\alpha}\varepsilon_k+O(\varepsilon^2_k).
\end{eqnarray*}
Using the previous equation and the fact that $\widetilde{f}(\gamma_k)=0$, we get
\begin{eqnarray*}
\widetilde{f}(\gamma_k)&=&\frac{4}{2-\alpha}\varepsilon_k
-\Frac{(\frac{1}{2}+\nu_{\alpha})(\frac{3}{2}+\nu_{\alpha})}{2i}\left( \frac{2}{2-\alpha}\imath\right)^{-1}
\frac{(-2)}{\left(-\frac{2-\alpha}{2}\imath k \pi \right)}\\\\
&-&\frac{1-\alpha}{i}
\Frac{(-2i)}{\left(-\frac{2-\alpha}{2}\imath k \pi \right)}
+O(\varepsilon_k^2)+ O\left(\frac{\varepsilon_k}{k}\right)+O\left(\frac{1}{k^{3-2\tilde{\alpha}}}\right)=0.\\ \\
\end{eqnarray*}
Hence
$$
\varepsilon_k=i\Frac{(2-\alpha)}{4}\left((\frac{1}{2}+\nu_{\alpha})(\frac{3}{2}+\nu_{\alpha})-\frac{4(1-\alpha)}{(2-\alpha)}\right)\Frac{1}{k\pi}
+ O(\varepsilon_k^2)+ O\left(\frac{\varepsilon_k}{k}\right)+O\left(\frac{1}{k^{3-2\tilde{\alpha}}}\right),
$$
We can write
\begin{equation}\label{332}
\gamma_k=-\frac{2-\alpha}{2}\imath\left(k+\frac{\nu_{\alpha}}{2}+\frac{5}{4}\right)\pi+
i\Frac{(2-\alpha)}{4}\left((\frac{1}{2}+\nu_{\alpha})(\frac{3}{2}+\nu_{\alpha})-\frac{4(1-\alpha)}{(2-\alpha)}\right)\Frac{1}{k\pi}
+\tilde{\varepsilon}_k,
\end{equation}
where $\tilde{\varepsilon}_k= o\left(\Frac{1}{k}\right)$.
Substituting (\ref{332}) into $\widetilde{f}(\gamma_k)=0$,
we get
$$
\widetilde{f}(\gamma_k)=\frac{4}{2-\alpha}\tilde{\varepsilon}_k+\rho i^{\tilde{\alpha}-1}\Frac{(-2i)}{\left(-\frac{2-\alpha}{2}\imath k \pi \right)^{3-2\tilde{\alpha}}}
+o\left(\Frac{1}{k^2}\right)+O(\tilde{\varepsilon}_k^2)=0.
$$
Hence
$$
\tilde{\varepsilon}_k=-\rho i^{1-\tilde{\alpha}}\left(\frac{2}{2-\alpha}\right)^{2-2\tilde{\alpha}}\Frac{1}{\left(k \pi \right)^{3-2\tilde{\alpha}}}
+o\left(\frac{1}{k^{3-2\tilde{\alpha}}}\right).
$$
It follows that
$$
\matrix{\gamma_k=-\frac{2-\alpha}{2}\imath\left(k+\frac{\nu_{\alpha}}{2}+\frac{5}{4}\right)\pi+
i\Frac{(2-\alpha)}{4}\left((\frac{1}{2}+\nu_{\alpha})(\frac{3}{2}+\nu_{\alpha})-\frac{4(1-\alpha)}{(2-\alpha)}\right)\Frac{1}{k\pi}\hfill &\cr
-\rho i^{1-\tilde{\alpha}}\left(\frac{2}{2-\alpha}\right)^{2-2\tilde{\alpha}}\Frac{1}{\left(k \pi \right)^{3-2\tilde{\alpha}}}
+o\left(\frac{1}{k^{3-2\tilde{\alpha}}}\right). &\cr}
$$
Since $\gamma_k^2=-\imath\lambda_k$, then
\begin{eqnarray*}
\lambda_k&=&\imath\gamma_k^2\\\\
&=&\imath\left[- C_0^2(k\pi)^2-C_1^2\pi^2-2 C_0 C_1 k\pi^2-2{C_0 C_2}
+2i^{2-\tilde{\alpha}}\Frac{C_0 C_3}{(k\pi)^{2-2\tilde{\alpha}}}
+o\left(\frac{1}{k^{2-2\tilde{\alpha}}}\right)\right]\\\\
&=&-i \left[C_0^2(k\pi)^2+C_1^2\pi^2+2 C_0 C_1 k\pi^2+2{C_0 C_2}\right]
-2i^{1-\tilde{\alpha}}\Frac{C_0 C_3}{(k\pi)^{2-2\tilde{\alpha}}}+o\left(\frac{1}{k^{2-2\tilde{\alpha}}}\right).
\end{eqnarray*}
where
$$
C_0=-\frac{2-\alpha}{2}, C_1=-\frac{2-\alpha}{2}\left(\frac{\nu_{\alpha}}{2}+\frac{5}{4}\right),
$$
$$
C_2=\Frac{(2-\alpha)}{4}\left((\frac{1}{2}+\nu_{\alpha})(\frac{3}{2}+\nu_{\alpha})-\frac{4(1-\alpha)}{(2-\alpha)}\right)
$$
and
$$
C_3=-\rho \left(\frac{2}{2-\alpha}\right)^{2-2\tilde{\alpha}}.
$$

\hfill$\diamondsuit$\\

$\bullet$\ Case 2: $\tilde{\alpha}< \frac{1}{2}$.
$$
\gamma_k=\gamma_k^0+\Frac{il}{k\pi}+\varepsilon_k,
$$
where $l=\Frac{(2-\alpha)}{4}\left((\frac{1}{2}+\nu_{\alpha})(\frac{3}{2}+\nu_{\alpha})-\frac{4(1-\alpha)}{(2-\alpha)}\right)$.
then
\begin{eqnarray*}
e^{2\imath(z-(\nu_{\alpha}+1)\frac{\pi}{2}-\frac{\pi}{4})}+1&=&1-e^{-\frac{4}{2-\alpha}\left(\varepsilon_k+\Frac{il}{k\pi}\right)}\\
&=&\frac{4}{2-\alpha}\left(\varepsilon_k+\Frac{il}{k\pi}\right)
-\Frac{1}{2}\left(\frac{4}{2-\alpha}\right)^2\Frac{(il)^2}{(k\pi)^2}
+O(\varepsilon^2_k)+O\left(\Frac{\varepsilon_k}{k}\right).
\end{eqnarray*}
\begin{eqnarray*}
e^{2\imath(z-(\nu_{\alpha}+1)\frac{\pi}{2}-\frac{\pi}{4})}-1&=&-2+\frac{4}{2-\alpha}\Frac{il}{k\pi}+O(\varepsilon_k).
\end{eqnarray*}
Using the previous equation and the fact that $\widetilde{f}(\gamma_k)=0$, we get
\begin{eqnarray*}
\widetilde{f}(\gamma_k)&=&\frac{4}{2-\alpha}\varepsilon_k+\frac{4}{2-\alpha}\left(\Frac{il}{k\pi}\right)
-\Frac{1}{2}\left(\frac{4}{2-\alpha}\right)^2\Frac{(il)^2}{(k\pi)^2}\\
&-&\Frac{(\frac{1}{2}+\nu_{\alpha})(\frac{3}{2}+\nu_{\alpha})}{2i}\left( \frac{2}{2-\alpha}\imath\right)^{-1}
\frac{(-2+\frac{4}{2-\alpha}\left(\Frac{il}{k\pi}\right)+O(\varepsilon_k))}{rk(1+\frac{\tilde{r}}{rk}+\frac{il}{rk^2\pi}+\frac{\varepsilon_k}{k})}\\ \\
&+&\frac{1-\alpha}{i}
\Frac{\left(2i+\frac{4}{2-\alpha}\Frac{l}{k\pi}+O(\varepsilon_k)\right)}{rk(1+\frac{\tilde{r}}{rk}+\frac{il}{rk^2\pi}+\frac{\varepsilon_k}{k})}
+O(\varepsilon_k^2)+ O\left(\frac{\varepsilon_k}{k}\right)+O\left(\frac{1}{k^{3-2\tilde{\alpha}}}\right)=0.\\ \\
&=&\frac{4}{2-\alpha}\varepsilon_k+\frac{4}{2-\alpha}\left(\Frac{il}{k\pi}\right)-\Frac{1}{2}\left(\frac{4}{2-\alpha}\right)^2\Frac{(il)^2}{(k\pi)^2}\\ \\
&-&\Frac{(\frac{1}{2}+\nu_{\alpha})(\frac{3}{2}+\nu_{\alpha})}{2i}\left( \frac{2}{2-\alpha}\imath\right)^{-1}
(-\Frac{2}{rk}+\Frac{2\tilde{r}}{r^2k^2}+\frac{4}{2-\alpha}\Frac{il}{rk^2\pi}+o(\frac{1}{k^2}))\\ \\
&&+(1-\alpha)(\Frac{2}{rk}-\Frac{2\tilde{r}}{r^2k^2}-\frac{4}{2-\alpha}\Frac{il}{rk^2\pi}+o(\frac{1}{k^2}))\\
\end{eqnarray*}
Hence
\begin{eqnarray*}
\frac{4}{2-\alpha}\varepsilon_k-\Frac{(\frac{1}{2}+\nu_{\alpha})(\frac{3}{2}+\nu_{\alpha})}{2i}\left( \frac{2}{2-\alpha}\imath\right)^{-1}\Frac{2\tilde{r}}{r^2k^2}
+(1-\alpha)\Frac{2\tilde{r}}{r^2k^2}+o(\frac{1}{k^2})=0.
\end{eqnarray*}
Then
$$
\varepsilon_k=\frac{2-\alpha}{4}(-\Frac{(\frac{1}{2}+\nu_{\alpha})(\frac{3}{2}+\nu_{\alpha})}{4}(2-\alpha)+(1-\alpha))\Frac{2\tilde{r}}{r^2k^2}+o(\frac{1}{k^2}).
$$
Now, we can write
\begin{equation}\label{89p}
\gamma_k=\gamma_k^0+\Frac{il}{k\pi}+\frac{2-\alpha}{4}\left(-\Frac{(\frac{1}{2}+\nu_{\alpha})(\frac{3}{2}+\nu_{\alpha})}{4}(2-\alpha)+(1-\alpha)\right)\Frac{2\tilde{r}}{r^2k^2}
+\varepsilon_k,
\end{equation}
Using the fact that $\widetilde{f}(\gamma_k)=0$, we get
$$
\frac{4}{2-\alpha}\varepsilon_k+\rho i^{\tilde{\alpha}-1}\Frac{(-2i)}{\left(-\frac{2-\alpha}{2}\imath k \pi \right)^{3-2\tilde{\alpha}}}+o\left(\frac{1}{k^{3-2\tilde{\alpha}}}\right)=0
$$
$$
\tilde{\varepsilon}_k=-\rho i^{1-\tilde{\alpha}}\left(\frac{2}{2-\alpha}\right)^{2-2\tilde{\alpha}}\Frac{1}{\left(k \pi \right)^{3-2\tilde{\alpha}}}
+o\left(\frac{1}{k^{3-2\tilde{\alpha}}}\right).
$$
Since $\gamma_k^2=-\imath\lambda_k$, then
\begin{eqnarray*}
\lambda_k&=&\imath\gamma_k^2\\\\
&=&\imath\left[- C_0^2(k\pi)^2-C_1^2\pi^2-\Frac {{C_2^2}}{(k\pi)^2}-2 C_0 C_1 k\pi^2-2{C_0 C_2}-2\Frac {{C_0 C_3}}{k}
-2\Frac {{C_1 C_2}}{k}\right.\\ \\
&&\left.+2i^{2-\tilde{\alpha}}\Frac{C_0 C_4}{(k\pi)^{2-2\tilde{\alpha}}}
+o\left(\frac{1}{k^{2-2\tilde{\alpha}}}\right)\right]\\ \\
&=&-i \left[C_0^2(k\pi)^2+C_1^2\pi^2+2 C_0 C_1 k\pi^2+2{C_0 C_2}+2\Frac {{C_0 C_3}}{k}
+2\Frac {{C_1 C_2}}{k}\right]\\ \\
&& -2i^{1-\tilde{\alpha}}\Frac{C_0 C_4}{(k\pi)^{2-2\tilde{\alpha}}}+o\left(\frac{1}{k^{2-2\tilde{\alpha}}}\right).
\end{eqnarray*}
where
$$
C_0=-\frac{2-\alpha}{2}, C_1=-\frac{2-\alpha}{2}\left(\frac{\nu_{\alpha}}{2}+\frac{5}{4}\right),
$$
$$
C_2=\Frac{(2-\alpha)}{4}\left((\frac{1}{2}+\nu_{\alpha})(\frac{3}{2}+\nu_{\alpha})-\frac{4(1-\alpha)}{(2-\alpha)}\right)
$$
$$
C_3=-\frac{2-\alpha}{4}\left(-\Frac{(\frac{1}{2}+\nu_{\alpha})(\frac{3}{2}+\nu_{\alpha})}{4}(2-\alpha)+(1-\alpha)\right)
\Frac{2C_1}{C_0^2}.
$$
and
$$
C_4=-\rho \left(\frac{2}{2-\alpha}\right)^{2-2\tilde{\alpha}}.
$$
Now, setting ${\tilde U}_k=(\lambda_k^0 I-{\cal A}) U_k$, where $U_k$ is
a normalized eigenfunction associated to $\lambda_k$ and
$$
\lambda_k^0=\left\{\matrix{-i \left[C_0^2(k\pi)^2+C_1^2\pi^2+2 C_0 C_1 k\pi^2+2{C_0 C_2}\right]\hfill & \hbox{ if }
\frac{1}{2} <\tilde{\alpha}<1,\hfill \cr
-i \left[C_0^2(k\pi)^2+C_1^2\pi^2+2 C_0 C_1 k\pi^2+2{C_0 C_2}+2\Frac {{C_0 C_3}}{k}+2\Frac {{C_1 C_2}}{k}\right]\hfill &
\hbox{ if } 0 <\tilde{\alpha}<\frac{1}{2}.\hfill \cr}\right.
$$
We then have
$$
\matrix{\|(\lambda_k^0 I-{\cal A})^{-1}\|_{{\cal L}({\cal H})}=\Sup_{U\in {\cal H}, U\not=0} \Frac{\|(\lambda_k^0 I-{\cal A})^{-1}U\|_{{\cal H}}}{\|U\|_{{\cal H}}}
&\geq &\Frac{\|(\lambda_k^0 I-{\cal A})^{-1}{\tilde U}_k\|_{{\cal H}}}{\|{\tilde U}_k\|_{{\cal H}}}\hfill \cr
&\geq &\Frac{\|U_k\|_{{\cal H}}}{\|(\lambda_k^0 I-{\cal A}) U_k\|_{{\cal H}}}.\hfill \cr}
$$
Hence, by Lemma \ref{l4.1}, we deduce that
$$
\|(\lambda_k^0 I-{\cal A})^{-1}\|_{{\cal L}({\cal H})}\geq
c |k|^{2-2\tilde{\alpha}}\equiv |\lambda_k^0|^{1-\tilde{\alpha}},
$$
Thus, the second condition of Lemma \ref{4.1} is not satisfied for $0< \tilde{\alpha}<1$. So that, the semigroup
$e^{t {\cal A}}$ is not exponentially stable. Thus the proof is complete.

\end{sloppypar}
\end{document}